# A Comparison Between Laguerre, Hermite, and Sinc Orthogonal Functions

Fattaneh Bayatbabolghani
Computer Science and Engineering
University of Notre Dame
fbayatba@nd.edu

Kourosh Parand
Computer Science
Shahid Beheshti University
k_parand@sbu.ac.ir

**Abstract**

A series of problems in different fields such as physics and chemistry are modeled by differential equations. Differential equations are divided into partial differential equations and ordinary differential equations which can be linear or nonlinear. One approach to solve those kinds of equations is using orthogonal functions into spectral methods. In this paper, we firstly describe Laguerre, Hermite, and Sinc orthogonal functions. Secondly, we select three interesting problems which are modeled as differential equations over the interval $[0, +\infty)$. Then, we use the collocation method as a spectral method for solving those selected problems and compare the performance of Laguerre, Hermite, and Sinc orthogonal functions in solving those types of equations.

## 1 Introduction

Mathematical modeling of many problems in science and engineering is in a form of differential equations. For some of those equations, there is not precise solution, and finding the best approximation for the solution is challenging. In this paper, we selected three problems which have no precise solution, and they are define in semi-infinite interval. The selected problems are as follows:

- Steady flow of a third grade fluid in a porous half space;
- Thomas-Fermi;
- The case of heat transfer equations arising in porous medium.

We propose solutions for the selected problems based on collocation method as a spectral method. The solutions are based on Modified generalized Laguerre, Hermite, and Sinc functions as orthogonal functions. We also, compare the proposed solutions to the existing solutions to show the accuracy of the proposed approaches.

For the rest of this paper, we review related literature in Section 2. Then, we explain all necessary background about collocation method and Modified generalized Laguerre, Hermite, and Sinc functions in Section 3. In Section 4, the details of the three selected applications are discussed, we explain our approached, and report the performance evaluation of the solutions. Finally, concluding remarks are presented in Section 5.

## 2  Related Work

Many of the current science and engineering problems are set in unbounded domains. In the context of spectral methods such as collocation and Galerkin methods [20], a number of approaches for treating unbounded domains have been proposed and investigated. The most common method is the use of polynomials that are orthogonal over unbounded domains, such as the Hermite and Laguerre spectral method [6, 7, 8, 10, 14, 21, 41, 42].

Guo [11, 66, 12, 13] proposed another method that proceeds by mapping the original problem in an unbounded domain to a problem in a bounded domain, and then using suitable Jacobi polynomials such as Gegenbauer polynomials to approximate the resulting problems. The Jacobi polynomials are a class of classical orthogonal polynomials and the Gegenbauer polynomials. The Legendre and Chebyshev polynomials are special cases of those orthogonal polynomials which have been used in several literature for solving differential equations [52, 19].

In [53], another method was introduced which is named domain truncation [53]. In domain truncation, infinite domain is replaced with $[-L, L]$ and semi-infinite interval with $[0, L]$ by choosing $L$ sufficiently large. There is another effective direct approach for solving such problems which is based on rational approximations. Christov [67] and Boyd [54, 55] developed spectral methods on unbounded intervals by using mutually orthogonal systems of rational functions. Boyd [54] defined a new spectral basis, named rational Chebyshev functions on the semi-infinite interval, by mapping to the Chebyshev polynomials. Guo et al. [15] introduced a new set of rational Legendre functions which are mutually orthogonal in $L^2(0, \infty)$. They applied a spectral scheme using the rational Legendre functions for solving the Korteweg-de Vries equation on the half-line. Boyd et al. [68] applied pseudo-spectral methods on a semi-infinite interval and compared rational Chebyshev, Laguerre and mapped Fourier sine methods. Parand et al. [30, 31, 32, 33, 34, 35, 63] also applied spectral method to solve nonlinear ordinary differential equations defined over the interval $I = [0, \infty)$. Their approach was based on rational tau and collocation methods.

## 3  Preliminaries

In collocation method, a function $f(x)$ defined over the interval $I = [0, \infty)$ can be expanded as

$$f(x) = \sum_{i=0}^{\infty} a_i \phi_i(x), \qquad (1)$$

where $\phi$s are orthogonal functions and

$$a_i = \frac{<f, \phi_i>_w}{<\phi_i, \phi_i>_w}.$$

If the infinite series in Equation (1) is truncated with $N$ terms, then it can be written as [35]

$$f(x) \simeq \sum_{i=0}^{N-1} a_i \phi_i(x) = A^T \phi(x),$$

with

$$A = [a_0, a_1, a_2, ..., a_{N-1}]^T,$$

$$\phi(x) = [\phi_0(x), \phi_1(x), ..., \phi_{N-1}(x)]^T.$$

In the following of this section, we introduce Laguerre, Hermite, and Sinc orthogonal functions, and describe how they can be used for function approximation over the interval $I = [0, \infty)$.

**Modified generalized Laguerre functions.** In here, we describe the basic notions and working tools concerning orthogonal modified generalized Laguerre functions. Orthogonal modified generalized Laguerre functions have been widely used for numerical solutions of differential equations on infinite intervals. $L_n^\alpha(x)$ (generalized Laguerre polynomial) is the $n$th eigen-function of the Sturm-Liouville problem [6, 16, 35]:

$$x \frac{d^2}{dx^2} L_n^\alpha(x) + (\alpha + 1 - x) \frac{d}{dx} L_n^\alpha(x) + n L_n^\alpha(x) = 0,$$
$$x \in I = [0, \infty), \qquad n = 0, 1, 2, ....$$

The generalized Laguerre in polynomial manner is defined with the following recurrence formula:

$$L_0^\alpha(x) = 1,$$
$$L_1^\alpha(x) = 1 + \alpha - x,$$
$$n L_n^\alpha(x) = (2n - 1 + \alpha - x) L_{n-1}^\alpha(x) - (n + \alpha - 1) L_{n-2}^\alpha(x).$$

These are orthogonal polynomials for the weight function $w_\alpha = x^\alpha e^{-x}$. We define Modified generalized Laguerre functions $\phi_j$ over the interval $(0, +\infty)$ as follows [35]:

$$\phi_j(x) = \exp(\frac{-x}{2L}) L_j^1(\frac{x}{L}), \quad L > 0.$$

This system is an orthogonal basis [9, 46] with weight function $w(x) = \frac{x}{L}$ and orthogonality property [35]:

$$< \phi_m, \phi_n >_{w_L} = (\frac{\Gamma(n + 2)}{L^2 n!}) \delta_{nm},$$

where $\delta_{nm}$ is the Kronecker function.
Laguerre-Gauss-Radau points and generalized Laguerre-Gauss-type interpolation were introduced by [17, 35, 47, 50]. Let

$$\Re_N = span\{1, x, ..., x^{2N-1}\},$$

we choose the collocation points relative to the zeroes of the functions [35]

$$p_j(x) = \phi_j(x) - (\frac{j+1}{j}) \phi_{j-1}(x).$$

Let $w(x) = \frac{x}{L}$ and $x_j$, $j = 0, 1, ..., N - 1$, be the $N$ modified generalized Laguerre functions-Radau points. The relation between modified generalized Laguerre functions orthogonal systems and modified generalized Laguerre functions integration is as follows [35, 45]:

$$\int_0^\infty f(x)w(x)\mathrm{d}x = \sum_{j=0}^{N-1} f_j(x)w_j + (\frac{\Gamma(N+2)}{(N)!(2N)!})f^{2N}(\xi)e^\xi,$$

where $0 < \xi < \infty$ and $w_j = x_j \frac{\Gamma(N+2)}{(L(N+1)![(N+1)\phi_{N+1}(x_j)]^2)}$, $j = 0, 1, 2, ..., N-1$. In particular, the second term on the right-hand side vanishes when $f(x)$ is a polynomial of degree at most $2N-1$ [35]. We define

$$I_N u(x) = \sum_{j=0}^{N-1} a_j \phi_j(x),$$

it as: $I_N u(x_j) = u(x_j)$, $j = 0, 1, 2, ..., N-1$. $I_N u$ is the orthogonal projection of u upon $\mathfrak{R}_N$ with respect to the discrete inner product and discrete norm as [35]:

$$<u, v>_{w,N} = \sum_{j=0}^{N-1} u(x_j)v(x_j)w_j,$$

$$\| u \|_{w,N} = <u, v>_{w,N}^{\frac{1}{2}},$$

thus for the modified generalized Laguerre functions Gauss-Radau interpolation we have

$$<I_N u, v>_{w,N} = <u, v>_{w,N}, \quad \forall u.v \in \mathfrak{R}_N.$$

**Hermite Functions.** Now, we elaborate the properties of Hermite functions. The reason of our focus on Hermite functions is that Hermite polynomials are generally not suitable in practice due to their wild asymptotic behavior at infinities [57]. Therefore, for the rest if this section, we only consider Hermite Functions.

The normalized Hermite functions of degree $n$ is defined by [58]

$$\widetilde{H}_n = \frac{1}{\sqrt{2^n n!}} e^{\frac{-x^2}{2}} H_n(x), \quad n \geq 0, x \in \mathbb{R}$$

where $\{\widetilde{H}_n\}$ is an orthogonal system in $L^2(\mathbb{R})$.

Also, the Hermite functions are well behaved with the decay property which is:

$$|\widetilde{H}_n(x)| \to 0, \quad as \ |x| \to \infty.$$

In addition, the three-term recurrence relation of Hermite functions implies [58]:

$$\widetilde{H}_{n+1}(x) = x\sqrt{\frac{2}{n+1}}\widetilde{H}_n(x) - \sqrt{\frac{n}{n+1}}\widetilde{H}_{n-1}(x), \quad n \geq 1,$$

$$\widetilde{H}_0(x) = e^{\frac{-x^2}{2}},$$

$$\widetilde{H}_1(x) = \sqrt{2}x e^{\frac{-x^2}{2}}.$$

More information about Hermite functions can be found in [58, 64, 65].

Note that, Hermite functions are defined over the interval $(-\infty, +\infty)$. For solving problems which are defined over the interval $(0, +\infty)$, we need to transform Hermite functions. One of the approaches to construct an approximation over the interval $(0, +\infty)$ is using mapping that is changing variable of the form [58]:

$$w = \Phi(z) = \frac{1}{k} ln(z),$$

where $k$ is a constant. Now, the transformed Hermite functions are

$$\widehat{H}_n(x) \equiv \widetilde{H}_n(x) o \Phi(x) = \widetilde{H}_n(\Phi(x)),$$

and the inverse map of $w = \Phi(z)$ is

$$z = \Phi^{-1}(w) = e^{kw}.$$

Therefor, we may define the inverse images of the spaced nodes $\{x_j\}_{x_j=-\infty}^{x_j=+\infty}$ as [58]:

$$\Gamma = \{\Phi^{-1}(t) : -\infty < t < +\infty\} = (0, +\infty),$$

and

$$\widetilde{x}_j = \Phi^{-1}(x_j) = e^{x_j}, \quad j = 0, 1, 2, \ldots$$

Let $w(x)$ denotes a non-negative, integrable, real-valued function over the interval $\Gamma$, We define

$$L_w^2(\Gamma) = \{v : \Gamma \to \mathbb{R} | v \text{ is measurable and } \|v\|_w < \infty\},$$

where

$$\|v\|_w = (\int_0^\infty |v(x)|^2 w(x) dx)^{\frac{1}{2}},$$

is the norm induced by the inner product of the space $L_w^2(\Gamma)$ [58],

$$< u, v >_w = \int_0^\infty u(x) v(x) w(x) dx.$$

Thus, $\{\widehat{H}_n(x)\}_{n \in \mathbb{N}}$ denotes a system which is mutually orthogonal

$$\langle \widehat{H}_n(x), \widehat{H}_m(x) \rangle_{w(x)} = \sqrt{\pi} \delta_{nm}.$$

This system is complete in $L_w^2(\Gamma)$. Therefore, for any function $f \in L_w^2(\Gamma)$, the following expansion holds [58]

$$f(x) \cong \sum_{k=-N}^{+N} f_k \widehat{H}_k(x),$$

with

$$f_k = \frac{\langle f(x), \widehat{H}_k(x)\rangle_{w(x)}}{\|\widehat{H}_k(x)\|^2_{w(x)}}$$

Now we define an orthogonal projection based on the transformed Hermite functions as given below [58]. Let set

$$\widehat{H}_N = span\{\widehat{H}_0(x), \widehat{H}_1(x), \ldots, \widehat{H}_n(x)\}.$$

The $L^2(\Gamma)$-orthogonal projection $\widehat{\xi}_N : L^2(\Gamma) \to \widehat{H}_N$ is a mapping in a way that for any $y \in L^2(\Gamma)$ [58],

$$\langle \widehat{\xi}_N y - y, \phi\rangle = 0 \quad \forall \phi \in \widehat{H}_N,$$

or equivalently,

$$\widehat{\xi}_N y(x) = \sum_{i=0}^{N} \widehat{a}_i \widehat{H}_i(x).$$

**Sinc Function.** Sinc functions are defined over interval $(-\infty, \infty)$ by [83] as follows:

$$Sinc(x) = \begin{cases} \frac{\sin(\pi x)}{\pi x} & x \neq 0 \\ 1 & x = 0 \end{cases}$$

For each integer $k$ and the mesh size $h$, the Sinc functions are defines on $\mathbb{R}$ by [84] as:

$$S_k(h, x) \equiv Sinc(\frac{x - kh}{h}) = \begin{cases} \frac{\sin(\frac{\pi}{h}(x - kh))}{\frac{\pi}{h}(x - kh)} & x \neq kh \\ 1 & x = kh \end{cases}$$

To solve a problem over the interval of $[0, \infty)$ by Sinc functions, we need to transform Sinc functions. Therefore, we can change a variable of the form:

$$w = \Phi(z) = ln(\sinh(z)).$$

The basis functions over $[0, \infty)$ are taken to be composite translates Sinc functions [83]:

$$S_k(x) \equiv S(k, h)(\Phi(x)) = Sinc(\frac{\Phi(x) - kh}{h}).$$

Then, the inverse map of $w = \Phi(z)$ is [85]

$$z = \Phi^{-1}(w) = ln(e^w + \sqrt{e^{2w} + 1}).$$

Thus,

$$x_k = \Phi^{-1}(kh) = ln(e^{kh} + \sqrt{e^{2kh} + 1}), \quad k = 0, \pm 1, \pm 2, \ldots.$$

Let $w(x)$ denotes a non-negative,integral-able, real-valued function over the interval $[0, \infty)$, we can define:

$$L^2_w(\Gamma) = \{v : \Gamma \to \mathbb{R} | v \text{ is measurable and } \|v\|_w < \infty\},$$

where

$$\|v\|_w = (\int_0^\infty |v(x)|^2 w(x) dx)^{\frac{1}{2}},$$

is the norm induced by the inner product of the space $L^2_w(\Gamma)$,

$$< u, v >_w = \int_0^\infty u(x) v(x) w(x) dx.$$

Thus, $\{S_k(x)\}_{k \in Z}$ denotes a system which is mutually orthogonal

$$\langle S_{k_n}(x), S_{k_m}(x) \rangle_{w(x)} = h S_{nm}.$$

Now, for any function $L^2_w(\Gamma)$, the following expansion holds [83]:

$$f(x) \cong \sum_{k=-N}^{+N} f_k S_k(x).$$

In addition, the $n$th derivation of the function $f$ at point $x_k$ can be approximated as:

$$\delta_{k,j}^{(0)} = [S(k,h)(\Phi(x))]|_{x=x_j} = \begin{cases} 0 & k \neq j \\ 1 & k = j \end{cases}$$

$$\delta_{k,j}^{(1)} = \frac{d}{d\Phi}[S(k,h)(\Phi(x))]|_{x=x_j} = \frac{1}{h} \begin{cases} \frac{(-1)^{j-k}}{j-k} & k \neq j \\ 0 & k = j \end{cases}$$

$$\delta_{k,j}^{(2)} = \frac{d^2}{d\Phi^2}[S(k,h)(\Phi(x))]|_{x=x_j} = \frac{1}{h^2} \begin{cases} \frac{-2(-1)^{j-k}}{(j-k)^2} & k \neq j \\ \frac{-\pi^2}{3} & k = j \end{cases}$$

# 4 Application of Laguerre, Hermite, and Sinc Orthogonal Functions in Solving Differential Equations

In this section, we describe three well-known problems which are modeled in the form of differential equations. Then, we solve each of them by Laguerre, Hermite, and Sinc Orthogonal Functions based on collocation methods [59, 60, 61, 62].

## 4.1 Steady Flow of a Third Grade Fluid in a Porous Half Space

In this section we focus on [51] which has discussed the flow of a third grade fluid in a porous half space. Based on [51] for unidirectional flow of a second grade fluid, the following formulation can be written:

$$(\nabla p)_x = -\frac{\mu\varphi}{k}(1 + \frac{\alpha_1}{\mu}\frac{\partial}{\partial t})u,$$

Also, for second grade fluid, the following formulation can be written:

$$(\nabla p)_x = -\frac{\varphi}{k}[\mu u + \alpha_1\frac{\partial u}{\partial t} + 2\beta_3(\frac{\partial u}{\partial y})^2 u],$$

where $\mu$ is the dynamic viscosity, $u$ is the denote the fluid velocity and $p$ is the pressure, $k$ and $\varphi$, respectively represent the permeability and porosity of the porous half space which occupies the region $y > 0$ and $\alpha_1$, $\beta_3$ are material constants. Non dimensional fluid velocity $f$ and the coordinate $z$ can be defined:

$$z = \frac{V_0}{\nu}y, \qquad f(z) = \frac{u}{V_0},$$
$$V_0 = u(0), \qquad \nu = \frac{\mu}{\rho},$$

where $\nu$ and $V_0$ represent the kinematic viscosity.
The boundary value problem modeling the steady state flow of a third grade fluid in a porous half space becomes [51]:

$$f''(z) + b_1 f'^2(z)f''(z) - b_2 f(z)f'^2(z) - b_3 f(z) = 0, \qquad (2)$$
$$f(0) = 1, \qquad f(\infty) = 0,$$

where $b_1$, $b_2$ and $b_3$ are defined as:

$$b_1 = \frac{6\beta_3 V_0^4}{\mu\nu^2},$$
$$b_2 = \frac{2\beta_3\varphi V_0^2}{k\mu},$$
$$b_3 = \frac{\varphi\nu^2}{kV_0^2}.$$

Note that the parameters are not independent, since

$$b_2 = \frac{b_1 b_3}{3}.$$

The homotopy analysis method for solution of Equation (2) found in [51]. Later Ahmad gave the asymptotic form of the solution and utilize this information to develop a series solution [56].

### 4.1.1 Solving the Problem with Modified Generalized Laguerre Functions

To apply modified generalized Laguerre collocation method to Equation (2) with its boundary conditions, firstly we expand $f(z)$ as follows:

$$I_N f(z) = \sum_{j=0}^{N-1} a_j \phi_j(z),$$

To find the unknown coefficients $a_j$'s, we substitute the truncated series $f(z)$ into Equation (2) and its boundary conditions. Also, we define Residual function as follows:

$$\begin{aligned} Res(z) = \sum_{j=0}^{N-1} a_j \phi_j''(z) + b_1 (\sum_{j=0}^{N-1} a_j \phi_j'(z))^2 \sum_{j=0}^{N-1} a_j \phi_j''(z) \\ - b_2 \sum_{j=0}^{N-1} a_j \phi_j(z) (\sum_{j=0}^{N-1} a_j \phi_j'(z))^2 - b_3 \sum_{j=0}^{N-1} a_j \phi_j(z), \end{aligned} \quad (3)$$

$$\sum_{j=0}^{N-1} a_j \phi_j(0) = 1, \quad (4)$$

$$\sum_{j=0}^{N-1} a_j \phi_j(\infty) = 0. \quad (5)$$

By applying $z$ in Equation (3) with the $N$ collocation points which are roots of functions $L_\alpha^N$, we have $N$ equations that generates a set of $N$ nonlinear equations; also, we have one boundary equation in Equation (4). Now, all of these equations can be solved by Newton method for the unknown coefficients. Note that Equation (5) is always true; therefore, we do not need to apply this boundary condition.

### 4.1.2 Solving the Problem with Hermite Functions

For solving Steady Flow problem, we use $\frac{1}{k} ln(z)$ for changing variable. Also, because of boundary conditions, we set following function:

$$p(z) = \frac{1}{1 + \lambda z + z^2},$$

and $\lambda$ is constant. Finally, we have

$$\widehat{\xi}_N f(z) = p(z) + \widehat{\xi}_N f(z).$$

that

$$\widehat{\xi}_N f(x) = \sum_{i=0}^{N} \widehat{a}_i \widehat{H}_i(z).$$

To find the unknown coefficients $a_i$'s, we substitute the truncated series $\widehat{\xi}_N f(z)$ into Equation (2). Also, we define Residual function of the form

$$\begin{aligned} Res(z) = (p''(z) + \widehat{\xi}_N f''(z)) + b_1 (p'(z) + \widehat{\xi}_N f'(z))^2 (p''(z) + \widehat{\xi}_N f''(z)) \\ - b_2 (p(z) + \widehat{\xi}_N f(z))(p'(z) + \widehat{\xi}_N f'(z))^2 - b_3 (p(z) + \widehat{\xi}_N f(z)) = 0. \end{aligned} \quad (6)$$

By applying $z$ in Equation (6) with the $N$ collocation points which are roots of transformed Hermite function, we have $N$ equations that generates a set of $N$ nonlinear equations. Now, all of these equations can be solved by Newton method for the unknown coefficients.

### 4.1.3 Solving the Problem with Sinc Functions

For solving Steady Flow problem with Sinc functions, we use $ln(\sinh(z))$ to change the variable. Also, because of the boundary conditions, we set the following function:

$$p(z) = \frac{1}{1 + \lambda z + z^2},$$

and $\lambda$ is constant. Finally, we have

$$f(z) \cong f_N(z) = p(z) + u_N(z),$$

that

$$u_N(z) = \sum_{k=-N}^{+N} c_k \frac{z S_k(z)}{z^2 + 1}.$$

The collocation points are

$$z_j = ln(e^{jh} + \sqrt{1 + e^{2jh}}), \quad j = -N, \ldots, +N,$$

and the derivations of $u_N(z)$ are

$$u_N(z_j) = \frac{c_j z_j}{z_j^2 + 1},$$

$$u'_N(z_j) = \sum_{k=-N}^{+N} c_k((\frac{1}{1 + z_j^2} - \frac{2z_j^2}{(1 + z_j^2)^2})\delta_{k,j}^{(0)} + (\frac{z_j \Phi'(z_j)}{1 + z_j^2})\delta_{k,j}^{(1)}),$$

$$u''_N(z_j) = \sum_{k=-N}^{+N} c_k((\frac{-6z_j}{(1+z_j^2)^2} + \frac{8z_j^3}{(1+z_j^2)^3})\delta_{k,j}^{(0)}$$
$$+(\frac{2\Phi'(z_j)}{1+z_j^2} - \frac{4z_j^2 \Phi'(z_j)}{(1+z_j^2)^2} + \frac{z_j \Phi'(z_j)}{1+z_j^2})\delta_{k,j}^{(1)} + (\frac{z_j(\Phi'(z_j))^2}{1+z_j^2})\delta_{k,j}^{(2)}).$$

### 4.1.4 Performance Evaluation

We present the results of our approximation for $N = 20$, $\alpha = 1$ and $L = 0.99$ in modified generalized Laguerre functions (MGLF), $N = 16$, $k = 1.2$, and $\lambda = 0.678301$ in Hermite functions (HF), and $N = 17$, $h = 1$, and $\lambda = 0.47$ in Sinc functions (SF) for solving this problem for some typical values of parameters, $b_1 = 0.6$, $b_2 = 0.1$ and $b_3 = 0.5$. In this problem the numerical solution of $f'(0)$ is important. Ahmad [56] obtained $f'(0)$ by the shooting method and founded, correct to six decimal positions, $f'(0) = -0.678301$. We compare the present methods with numerical solution and Ahmad solution [56] in Table 1. It shows the present methods are highly accurate. Also, the solutions are presented graphically in Figures 1–3.

Table 1: Comparison between MGLF, HF, SH and Ahmad solutions [56] for $b_1 = 0.6$, $b_2 = 0.1$ and $b_3 = 0.5$.

| $z$ | Ahmad method[56] | MGLF | HF | SF | numerical[56] |
|---|---|---|---|---|---|
| 0.0 | 1.00000 | 1.00000 | 1.00000 | 1.00000 | 1.00000 |
| 0.2 | 0.87220 | 0.87261 | 0.87261 | 0.87278 | 0.87260 |
| 0.4 | 0.76010 | 0.76063 | 0.76064 | 0.76035 | 0.76060 |
| 0.6 | 0.66190 | 0.66243 | 0.66243 | 0.66178 | 0.66240 |
| 0.8 | 0.57600 | 0.57650 | 0.57647 | 0.57597 | 0.57650 |
| 1.0 | 0.50100 | 0.50144 | 0.50139 | 0.50115 | 0.50140 |
| 1.2 | 0.43560 | 0.43595 | 0.43591 | 0.43583 | 0.43590 |
| 1.6 | 0.32890 | 0.32920 | 0.32917 | 0.32905 | 0.32920 |
| 2.0 | 0.24820 | 0.24838 | 0.24839 | 0.24802 | 0.24840 |
| 2.5 | 0.17440 | 0.17455 | 0.17459 | 0.17426 | 0.17450 |
| 2.7 | 0.15140 | 0.15156 | 0.15161 | 0.15141 | 0.15160 |
| 3.0 | 0.12250 | 0.12261 | 0.12270 | 0.12265 | 0.12260 |
| 3.6 | 0.08016 | 0.08024 | 0.08036 | 0.08025 | 0.08024 |
| 4.0 | 0.06042 | 0.06047 | 0.06060 | 0.06033 | 0.06047 |
| 4.2 | 0.05245 | 0.05250 | 0.05261 | 0.05233 | 0.05250 |
| 4.4 | 0.04553 | 0.04558 | 0.04567 | 0.04543 | 0.04558 |
| 4.6 | 0.03953 | 0.03957 | 0.03964 | 0.03948 | 0.03957 |
| 4.8 | 0.03432 | 0.03435 | 0.03440 | 0.03434 | 0.03435 |
| 5.0 | 0.02979 | 0.02982 | 0.02984 | 0.02987 | 0.02982 |
| $f'(0)$ | $-0.681835$ | $-0.678297$ | $-0.678301$ | $-0.677843$ | $-0.678301$ |

Figure 1: Graph of numerical approximation of $f(z)$ by MGLF.

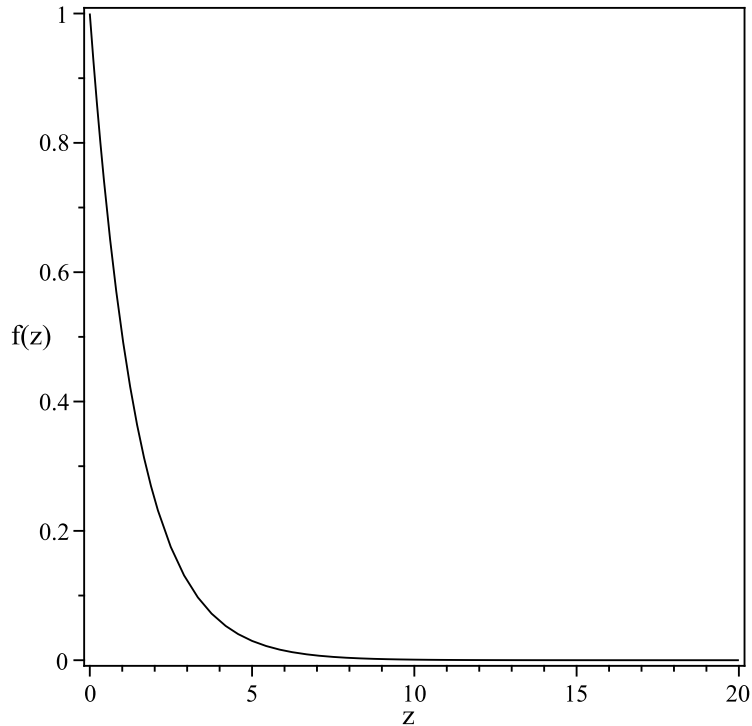

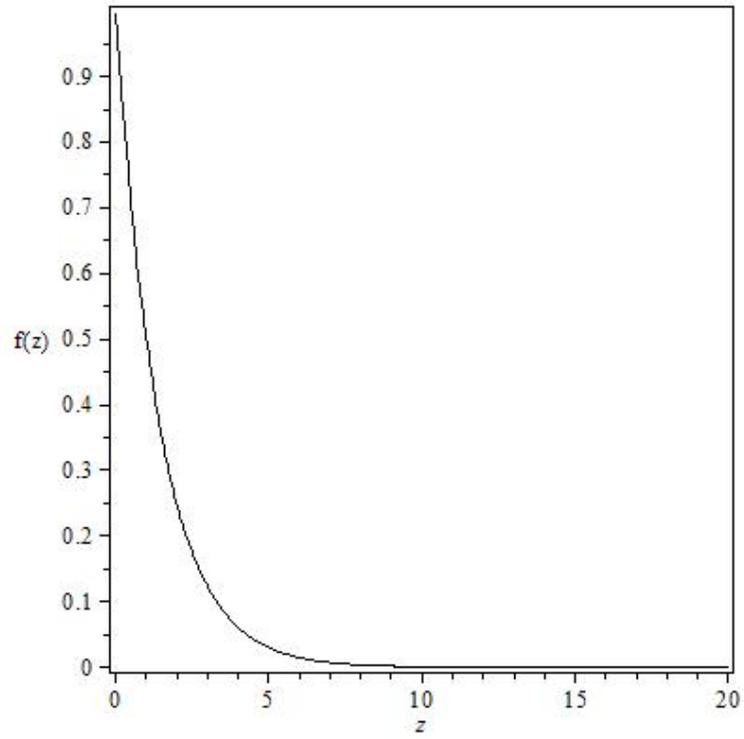

Figure 2: Graph of numerical approximation of $f(z)$ by HF.

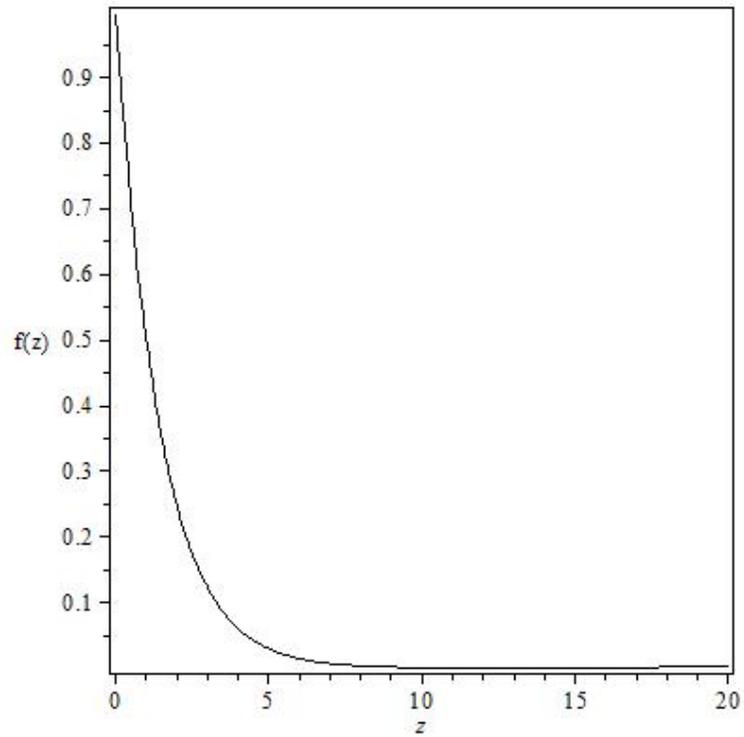

Figure 3: Graph of numerical approximation of $f(z)$ by SF.

## 4.2 Thomas-Fermi

Thomas-Fermi equation is One of the most important nonlinear ordinary differential equations that occurs in semi-infinite interval, as following [33, 69, 70]:

$$\frac{d^2y}{dx^2} = \frac{1}{\sqrt{x}}y^{\frac{3}{2}}(x), \tag{7}$$

which appears in the problem of determining the effective nuclear charge in heavy atoms. Also, it has following Boundary conditions:

$$y(0) = 1, \quad y(\infty) = 0.$$

The Thomas-Fermi equation is useful for calculating form factors and for obtaining effective potentials which can be used as initial trial potentials in self-consistent field calculations. The problem has been solved by different techniques [33, 71, 72, 73, 74, 75, 76, 77, 78, 79, 80]. [72, 73, 74] used perturbative approach to determine analytic solutions for Thomas-Fermi equation. Adomian [75] applied the decomposition method for solving Thomas-Fermi equation and then Wazwaz [76] proposed a non-perturbative approximate solution to this equation by using the modified decomposition method in a direct manner without any need to a perturbative expansion or restrictive assumptions. Liao [77] solved Thomas-Fermi equation by homotopy analysis method. Khan [78], used the homotopy analysis method (HAM) with a new and better transformation which improved the results in comparison with Liao's work. In [79], the quasilinearization approach was applied for solving Equation (7). This method approximated the solution of a nonlinear differential equation by treating the nonlinear terms as a perturbation about the linear ones, and unlike perturbation theories is not based on the existence of some kind of a small parameter. Ramos [80] presented two piecewise quasilinearization methods for the numerical solution of Equation (7). Both methods were based on the piecewise linearization of ordinary differential equations [33]. In addition, Parand [33] Solved Thomas-Fermi equation by Rational Chebyshev pseudo-spectral approach.

### 4.2.1 Solving the Problem with Modified Generalized Laguerre Functions

To apply modified generalized Laguerre collocation method to Equation (7) with its boundary conditions, firstly we expand $f(x)$ as follows:

$$I_N f(x) = \sum_{j=0}^{N-1} a_j \phi_j(x),$$

To find the unknown coefficients $a_j$'s, we substitute the truncated series $f(x)$ into Equation (7) and its boundary conditions. Also, we define Residual function as follows:

$$Res(x) = \sum_{j=0}^{N-1} a_j \phi_j''(x) - x^{\frac{1}{2}} (\sum_{j=0}^{N-1} a_j \phi_j'(x))^{\frac{3}{2}} \tag{8}$$

$$\sum_{j=0}^{N-1} a_j \phi_j(0) = 1, \tag{9}$$

$$\sum_{j=0}^{N-1} a_j \phi_j(\infty) = 0. \tag{10}$$

By applying $x$ in Equation (8) with the $N$ collocation points which are roots of functions $L_\alpha^N$, we have $N$ equations that generates a set of $N$ nonlinear equations; also, we have one boundary equation in Equation (9). Now, all of these equations can be solved by Newton method for the unknown coefficients. Note that Equation (10) is always true; therefore, we do not need to apply this boundary condition.

### 4.2.2 Solving the Problem with Hermite Functions

For solving Thomas-Fermi, we use $\frac{1}{k}ln(x)$ for changing variable. Also, because of boundary conditions, we set following function:

$$p(x) = \frac{1}{1 + \lambda x + x^2},$$

and $\lambda$ is constant. Finally, we have

$$\widehat{\xi}_N f(x) = p(x) + \widehat{\xi}_N f(x).$$

that

$$\widehat{\xi}_N f(x) = \sum_{i=0}^{N} \widehat{a}_i \widehat{H}_i(x).$$

To find the unknown coefficients $a_i$'s, we substitute the truncated series $\widehat{\xi}_N f(x)$ into Equation (7). Also, we define Residual function of the form

$$Res(x) = p''(x) + \widehat{\xi}_N f''(x) - x^{-\frac{1}{2}} \{p(x) + \widehat{\xi}_N f(x)\}^{\frac{3}{2}} = 0. \tag{11}$$

By applying $x$ in Equation (11) with the $N$ collocation points which are roots of transformed Hermite function, we have $N$ equations that generates a set of $N$ nonlinear equations. Now, all of these equations can be solved by Newton method for the unknown coefficients.

### 4.2.3 Solving the Problem with Sinc Functions

For solving Thomas-Fermi problem with Sinc functions, we use $ln(\sinh(x))$ to change the variable. Also, because of the boundary conditions, we set the following function:

$$p(x) = \frac{\lambda}{\lambda + x},$$

and $\lambda$ is constant. Finally, we can approximate $f(x)$ as:

$$f(x) \cong f_N(x) = p(x) + u_N(x),$$

that

$$u_N(x) = \sum_{k=-N}^{+N} c_k \frac{xS_k(x)}{x^2 + 1}.$$

The collocation points are

$$x_j = ln(e^{jh} + \sqrt{1+e^{2jh}}), \quad j = -N, \ldots, +N,$$

and the derivations of $u_N(x)$ are

$$u_N(x_j) = \frac{c_j x_j}{x_j^2 + 1},$$

$$u'_N(x_j) = \sum_{k=-N}^{+N} c_k\left(\left(\frac{1}{1+x_j^2} - \frac{2x_j^2}{(1+x_j^2)^2}\right)\delta_{k,j}^{(0)} + \left(\frac{x_j \Phi'(x_j)}{1+x_j^2}\right)\delta_{k,j}^{(1)}\right),$$

$$u'_N(x_j) = \sum_{k=-N}^{+N} c_k\left(\left(\frac{-6x_j}{(1+x_j^2)^2} + \frac{8x_j^3}{(1+x_j^2)^3}\right)\delta_{k,j}^{(0)}\right.$$
$$+ \left(\frac{2\Phi'(x_j)}{1+x_j^2} - \frac{4x_j^2 \Phi'(x_j)}{(1+x_j^2)^2} + \frac{x_j \Phi'(x_j)}{1+x_j^2}\right)\delta_{k,j}^{(1)} + \left(\frac{x_j(\Phi'(x_j))^2}{1+x_j^2}\right)\delta_{k,j}^{(2)}\right).$$

### 4.2.4 Performance Evaluation

The initial slope $y'(0)$ of the Thomas-Fermi equation is calculated by Kobayashi [81] as $y'(0) = -1.588071$. Table 2 shows the approximations of $y(x)$ and $y'(0)$ obtained by the present methods for $N = 7$, $\alpha = 1$ and $L = 0.675$ in modified generalized Laguerre functions (MGLF), $N = 15$, $k = 0.9$ and $\lambda = 1.588071$ in Hermite functions (HF), and $N = 11$, $h = 1$, and $\lambda = 0.77$ in Sinc functions (SF) to those obtained by Liao [82] and Kobayashi [81]. Figures 5–7 show the resulting graph of Thomas-Fermi which tends to zero as x increases by boundary condition $y(\infty) = 0$. They are compared with Liao's results that shown by square.

## 4.3 The Case of Heat Transfer Equations Arising in Porous Medium

Natural convective heat transfer in porous media has received considerable attention during the past few decades. This interest can be attributed due to its wide range of applications in ceramic processing, nuclear reactor cooling system, crude oil drilling, chemical reactor design, ground water pollution and filtration processes. External natural convection in a porous medium adjacent to heated bodies was analyzed by Nield and Bejan [28], Merkin [22, 23], Minkowycz and Cheng [24, 25, 26], Pop and Cheng [5, 36], Ingham and Pop [18]. In all of these analysis, it is assumed that boundary layer approximations are applicable and the coupled set of governing equations are solved by numerical methods. Also, [1, 40] worked on this problem. Parand [29] Compared two common collocation approaches based on radial basis functions for the case of heat transfer equations arising in porous medium.

In present study, we consider the problem of natural convection about an inverted heated cone embedded in a porous medium of infinite extent. No similarity solution exists for the truncated cone, but for the case of full cone similarity solutions exist if the prescribed wall temperature or surface heat flux is a power function of distance from the vertex of the inverted cone [5, 28, 43, 44]. Bejan and Khair [2] used Darcy's law to study the vertical natural convective flows driven by temperature and concentration gradients. Nakayama and Hossain [27] applied the integral method to obtain the heat and mass transfer by free convection from a vertical surface with constant wall temperature and concentration. Yih [48] examined the coupled heat and mass transfer by free convection over a truncated cone in porous media for variable wall temperature and variable heat and mass fluxes, Also he [49] applied the uniform transpiration effect on coupled heat and mass transfer in mixed convection about inclined surfaces in porous media for the entire regime. Cheng [3] uses an integral approach to study the heat and mass transfer by natural convection from truncated cones in porous media with variable wall temperature and [4] studies the Soret and Dufour effects on the boundary layer flow due to natural convection heat and mass transfer over a vertical cone in a porous medium saturated with Newtonian fluids with constant wall temperature. Natural convective mass transfer from upward-pointing vertical cones, embedded in saturated porous media, has been studied using the limiting diffusion [39]. The natural convection along an isothermal wavy cone embedded in a fluid-saturated porous medium are presented in [37, 38]. In [43, 44] fluid flow and heat transfer of vertical full cone embedded in porous media have been solved by Homotopy analysis method.

Table 2: Comparison between MGLF, HF, SH and Liao [82].

| $z$ | Liao [82] | MGLF | HF | SF |
|---|---|---|---|---|
| 00.25 | 0.755202000 | 0.765698401 | 0.754795330 | 0.755501513 |
| 00.50 | 0.606987000 | 0.611094841 | 0.606658908 | 0.606591374 |
| 00.75 | 0.502347000 | 0.504222143 | 0.502110510 | 0.502017836 |
| 01.00 | 0.424008000 | 0.426286491 | 0.423811203 | 0.424181905 |
| 01.25 | 0.363202000 | 0.366441192 | 0.363027725 | 0.363696263 |
| 02.00 | 0.243009000 | 0.245827159 | 0.242918233 | 0.243257878 |
| 02.25 | 0.215895000 | 0.217768255 | 0.215819818 | 0.216156024 |
| 02.50 | 0.192984000 | 0.193946429 | 0.192917948 | 0.193395401 |
| 02.75 | 0.173441000 | 0.173698712 | 0.173379623 | 0.174078699 |
| 03.00 | 0.156633000 | 0.156469166 | 0.156573773 | 0.157498937 |
| 03.25 | 0.142070000 | 0.141769745 | 0.142013368 | 0.143125471 |
| 03.50 | 0.129370000 | 0.129167788 | 0.129316613 | 0.130577926 |
| 03.75 | 0.118229000 | 0.118284289 | 0.118180209 | 0.119583902 |
| 04.00 | 0.108404000 | 0.108794792 | 0.108360441 | 0.109933372 |
| 08.00 | 0.036587300 | 0.035764064 | 0.036580427 | 0.046325607 |
| 15.00 | 0.010805400 | 0.009355939 | 0.010803774 | 0.049089345 |
| 20.00 | 0.005784940 | 0.001828955 | 0.005792831 | 0.036992990 |
| 30.00 | 0.005784940 | 0.000018510 | 0.002252634 | 0.024996549 |
| $y'$ | Kobayashi [81] | MGLF | HF | SF |
| $y'(0)$ | $-1.588071$ | $-1.158425$ | $-1.588071$ | $-1.580380$ |

Figure 4: Graph of numerical approximation of $f(x)$ by MGLF compare to Liao [82].

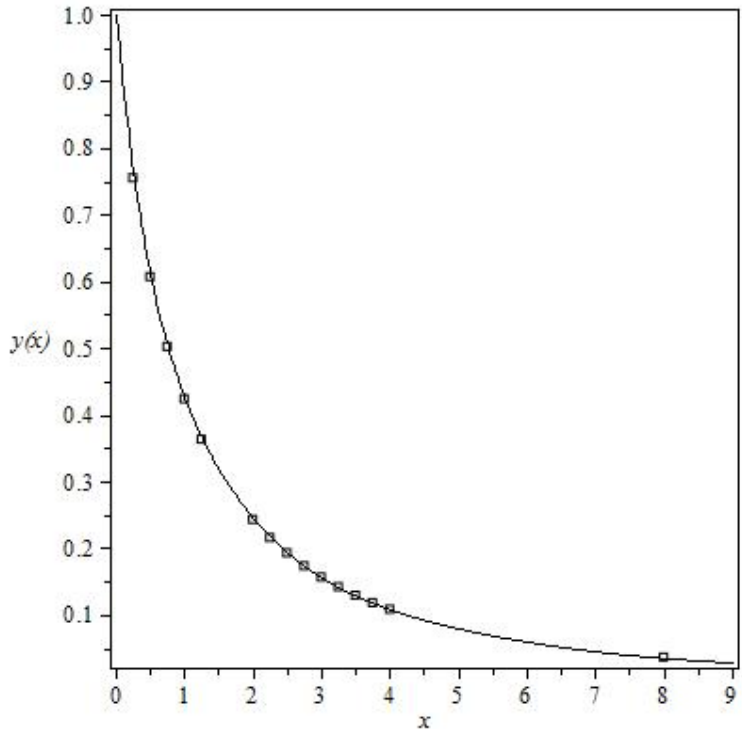

Figure 5: Graph of numerical approximation of $f(x)$ by HF compare to Liao [82].

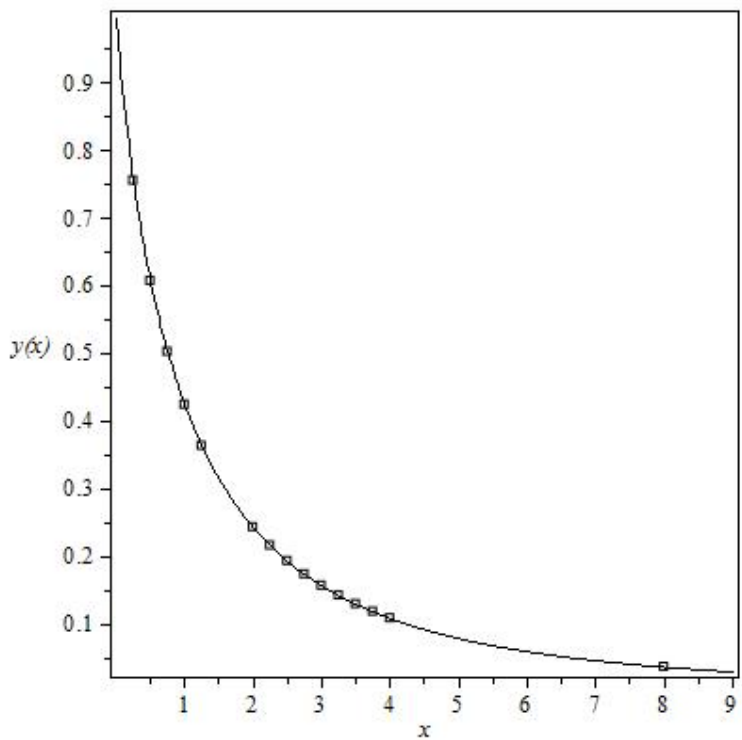

Figure 6: Graph of numerical approximation of $f(x)$ by SF compare to Liao [82].

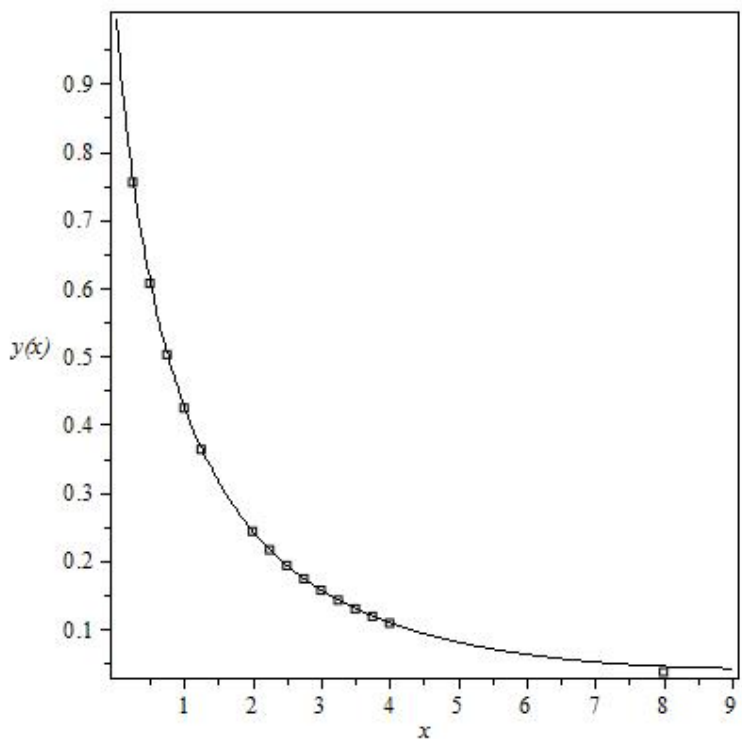

If we want to express the problem formulation of this model, we can consider an inverted cone with semi-angle $\gamma$ and take axes in the manner indicated in Figure 7(a). The boundary layer develops over the heated frustum $x = x_0$. In terms of the stream function $\psi$ defined by:

$$u = \frac{1}{r}\frac{\partial \psi}{\partial y}, \quad v = -\frac{1}{r}\frac{\partial \psi}{\partial x}.$$

The boundary layer equations for natural convection of Darcian fluid about a cone are:

$$\frac{\partial}{\partial x}(ru) + \frac{\partial}{\partial y}(rv) = 0$$

$$u = \frac{\rho_\infty \beta K g \cos\gamma (T - T_\infty)}{\mu}$$

$$\frac{1}{r}(\frac{\partial \psi}{\partial y}\frac{\partial T}{\partial x} - \frac{\partial \psi}{\partial x}\frac{\partial T}{\partial y}) = \alpha \frac{\partial^2 T}{\partial y^2}.$$

For a thin boundary layer, $r$ is obtained approximately $x\sin(\gamma)$. Suppose that a power law of heat flux is prescribed on the frustum. Accordingly, the boundary conditions at infinity are:

$$u = 0, \quad T = T_\infty, \quad y \to \infty,$$

and at the wall are

$$v = 0, \quad y = 0.$$

The surface heat flux $q_w$ is prescribed by

$$q_w = -k(\frac{\partial T}{\partial y})_{y=0} = A(x - x_0)^\lambda \quad x_0 \leq x \leq \infty.$$

For the case of a full cone ($x_0 = 0$, Figure 7(b)) a similarity solution exists. In the case of prescribed surface heat flux, we have:

$$\psi = \alpha r (Ra_x)^{1/3} f(\eta), \qquad (12)$$

$$T - T_\infty = \frac{q_w x}{k} (Ra_x)^{-\frac{1}{3}} \theta(\eta),$$

$$\eta = \frac{y}{x} (Ra_x)^{1/3},$$

where

$$Ra_x = \frac{\rho_\infty \beta g K \cos(\gamma) q_w x^2}{\mu \alpha k}, \qquad (13)$$

is the local Rayleigh number for the case of prescribed surface heat flux. The governing equations become

$$f' = \theta \qquad (14)$$

$$\theta'' + \frac{\lambda + 5}{2} f \theta' - \frac{2\lambda + 1}{3} f' \theta = 0.$$

Subjected to boundary conditions as:

$$f(0) = 0, \quad \theta'(0) = -1, \quad \theta(\infty) = 0. \qquad (15)$$

Finally from Equations (14) and (15) we have:

$$f''' + \left(\frac{\lambda + 5}{2}\right) f f'' - \left(\frac{2\lambda + 1}{3}\right) (f')^2 = 0, \qquad (16)$$

$$f(0) = 0, \quad f''(0) = -1, \quad f'(\infty) = 0.$$

It is of interest to obtain the value of the local Nusselt number which is defined as:

$$Nu_x = \frac{q_w x}{k(T_w - T_\infty)}. \qquad (17)$$

From Equations (17), (12) and (13) it follows that the local Nusselt number which is interest to obtain given by:

$$Nu_x = Ra_x^{1/3}[-\theta(0)].$$

### 4.3.1 Solving the Problem with Modified Generalized Laguerre Functions

To apply modified generalized Laguerre collocation method to Equation (16), at first we expand $f(\eta)$ as follows:

$$I_N f(\eta) = \sum_{j=0}^{N-1} a_j \phi_j(\eta),$$

Figure 7: (a) Coordinate system for the boundary layer on a heated frustum of a cone, (b) full cone, $x_0 = 0$.

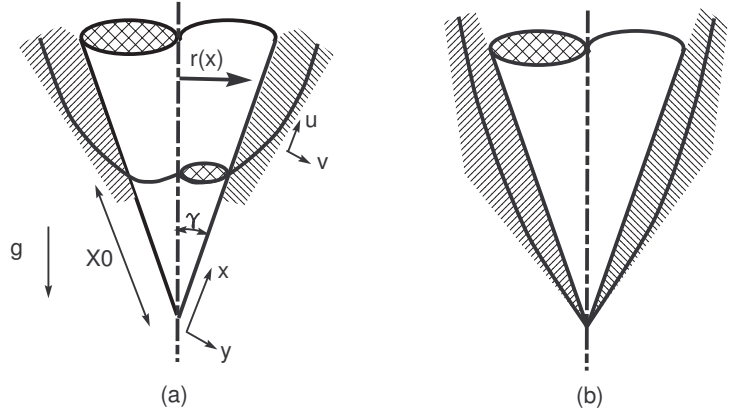

to find the unknown coefficients $a_j$'s, we substitute the truncated series $f(\eta)$ into Equation (16) and boundary conditions in it. Also, we define Residual function of the form

$$Res(\eta) = \sum_{j=0}^{N-1} a_j \phi_j'''(\eta) + (\tfrac{\lambda+5}{2}) \sum_{j=0}^{N-1} a_j \phi_j(\eta) \sum_{j=0}^{N-1} a_j \phi_j''(\eta) \\ - (\tfrac{2\lambda+2}{3})(\sum_{j=0}^{N-1} a_j \phi_j'(\eta))^2, \tag{18}$$

$$\sum_{j=0}^{N-1} a_j \phi_j(0) = 0, \quad \sum_{j=0}^{N-1} a_j \phi_j''(0) = -1, \tag{19}$$

$$\sum_{j=0}^{N-1} a_j \phi_j'(\infty) = 0. \tag{20}$$

By applying $\eta$ in Equation (18) with the $N$ collocation points which are roots of functions $L_\alpha^N$, we have $N$ equations that generates a set of $N$ nonlinear equations; also, we have two boundary equations in Equation (19). Now, all of these equations can be solved by Newton method for the unknown coefficients. We must mention Equation (20) is always true; therefore, we do not need to apply this boundary condition.

### 4.3.2 Solving the Problem with Hermite Functions

For solving this problem, we used $\frac{1}{k} ln(\eta)$ for changing variable. Also, because of boundary conditions, we set following function:

$$p(\eta) = \frac{\beta^2 \eta x}{2(\beta + \eta)},$$

and $\beta$ is constant. Finally, we have

$$\widehat{\xi}_N f(\eta) = p(\eta) + \widehat{\xi}_N f(\eta).$$

that

$$\widehat{\xi}_N f(\eta) = \sum_{i=0}^{N} \widehat{a}_i \widehat{H}_i(\eta).$$

To find the unknown coefficients $a_i$'s, we substitute the truncated series $\widehat{\xi}_N f(x)$ into Equation (16). Also, we define Residual function of the form

$$Res(\eta) = p'''(\eta) + \widehat{\xi}_N f'''(\eta) + \tfrac{\lambda+5}{2}(p(\eta) + \widehat{\xi}_N f(\eta))(p''(\eta) \qquad (21)$$
$$+\widehat{\xi}_N f''(\eta)) - \tfrac{2\lambda+1}{3}(p'(\eta) + \widehat{\xi}_N f'(\eta))^2.$$

By applying $\eta$ in Equation (21) with the $N$ collocation points which are roots of transformed Hermite function, we have $N$ equations that generates a set of $N$ nonlinear equations. Now, all of these equations can be solved by Newton method for the unknown coefficients.

### 4.3.3 Solving the Problem with Sinc Functions

For solving Thomas-Fermi problem with Sinc functions, we use $ln(\eta)$ to change the variable. Also, because of the boundary conditions, we set the following function:

$$p(\eta) = \frac{\beta^2 \eta}{2(\beta + \eta)},$$

and $\beta$ is constant. Finally, we can approximate $f(\eta)$ as:

$$f(\eta) \cong f_N(\eta) = p(\eta) + u_N(\eta),$$

that

$$u_N(\eta) = \sum_{k=-N}^{+N} c_k \frac{\eta^3 S_k(\eta)}{\eta^3 + 1}.$$

The collocation points are

$$\eta_j = e^{jh}, \quad j = -N, \ldots, +N,$$

and the derivations of $u_N(x)$ are

$$u_N(\eta_j) = \frac{c_j \eta_j^3}{\eta_j^3 + 1},$$

$$u'_N(\eta_j) = \sum_{k=-N}^{+N} c_k((\frac{3\eta_j^2}{1+\eta_j^3} - \frac{3\eta_j^5}{(1+\eta_j^3)^2})\delta_{k,j}^{(0)} + (\frac{\eta_j^3 \Phi'(\eta_j)}{1+\eta_j^3})\delta_{k,j}^{(1)}),$$

$$u''_N(\eta_j) = \sum_{k=-N}^{N} c_k((\frac{6\eta_j}{1+\eta_j^3} - \frac{24\eta_j^4}{(1+\eta_j^3)^2} + \frac{18\eta_j^7}{(1+\eta_j^3)^3})\delta_{k,j}^{(0)}$$
$$+(\frac{6\eta_j^2 \phi'(\eta_j)}{1+\eta_j^3} - \frac{6\eta_j^5 \phi'(\eta_j)}{(1+\eta_j^3)^2} + \frac{\eta_j^3 \phi''(\eta_j)}{1+\eta_j^3})\delta_{k,j}^{(1)}$$
$$+(\frac{\eta_j^3 (\phi'(\eta_j))^2}{1+\eta_j^3})\delta_{k,j}^{(2)}),$$

Table 3: Comparison of $f'(0)$ for various $\lambda$ between MGLF collocation method and Runge-Kutta solution.

| $\lambda$ | Runge-Kutta solution [44] | MGLF | $\alpha$ | $L$ |
|---|---|---|---|---|
| 0 | 0.94760 | 0.947697256 | 1 | 1.2985 |
| $\frac{1}{4}$ | 0.91130 | 0.911292295 | 0.94869 | 1.24093 |
| $\frac{1}{3}$ | 0.90030 | 0.900305806 | 1 | 1.15 |
| $\frac{1}{2}$ | 0.87980 | 0.879332090 | 1 | 1.09 |
| $\frac{3}{4}$ | 0.85220 | 0.852287074 | 0.04 | 1.0394 |
| 1 | 0.82760 | 0.825288795 | 0.655 | 1.115 |

Table 4: Comparison of $f'(0)$ for various $\lambda$ between HF collocation method and Runge-Kutta solution.

| $\lambda$ | Runge-Kutta solution [44] | HF | $k$ | $\beta$ |
|---|---|---|---|---|
| 0 | 0.94760 | 0.947350000 | 0.00005 | 1.8947 |
| $\frac{1}{4}$ | 0.91130 | 0.911300000 | 0.00005 | 1.8226 |
| $\frac{1}{3}$ | 0.90030 | 0.900350000 | 0.00005 | 1.8007 |
| $\frac{1}{2}$ | 0.87980 | 0.879330000 | 0.00005 | 1.75866 |
| $\frac{3}{4}$ | 0.85220 | 0.852100000 | 0.00005 | 1.7042 |
| 1 | 0.82760 | 0.827600000 | 0.00005 | 1.6552 |

$$\begin{aligned}
u'''_N(\eta_j) = \sum_{k=-N}^{N} c_k & \Big( \big( \tfrac{6}{1+\eta_j^3} - \tfrac{114\eta_j^3}{(1+\eta_j^3)^2} + \tfrac{270\eta_j^6}{(1+\eta_j^3)^3} - \tfrac{169x^9}{(1+\eta_j^3)^4} \big) \delta_{k,j}^{(0)} \\
& + \big( \tfrac{18\eta_j \phi'(\eta_j)}{1+\eta_j^3} - \tfrac{72\eta_j^4 \phi'(\eta_j)}{(1+\eta_j^3)^2} + \tfrac{54\eta_j^7 \phi'(\eta_j)}{(1+\eta_j^3)^3} + \tfrac{9\eta_j^2 \phi''(\eta_j)}{1+\eta_j^3} - \tfrac{9\eta_j^5 \phi''(\eta_j)}{(1+\eta_j^3)^2} + \tfrac{\eta_j^3 \phi'''(\eta_j)}{1+\eta_j^3} \big) \delta_{k,j}^{(1)} \\
& + \big( \tfrac{6\eta_j^2 \phi'(\eta_j)}{1+\eta_j^3} - \tfrac{6\eta_j^5 \phi'(\eta_j)}{(1+\eta_j^3)^2} + \tfrac{3\eta_j^3 \phi''(\eta_j)\phi'(\eta_j)}{1+\eta_j^3} + \tfrac{3\eta_j^2 (\phi'(\eta_j))^2}{1+\eta_j^3} - \tfrac{3\eta_j^5 (\phi'(\eta_j))^2}{(1+\eta_j^3)^2} \big) \delta_{k,j}^{(2)} \\
& + \big( \tfrac{\eta_j^3 (\phi'(\eta_j))^3}{1+\eta_j^3} \big) \delta_{k,j}^{(3)} \Big).
\end{aligned}$$

### 4.3.4 Performance Evaluation

In the following tables and figures we make a comparison between Runge-Kutta solution is obtained by the MATLAB software command ODE45 [44] with the obtained results of three presented method in this paper. The results for $f'(0)$ with the presented methods have been shown in Tables 3–5 and comparison has been made between the Runge-Kutta solution. Some of the computed results for the variations with $\eta$ of the functions $f'$ for $\lambda = 1/4$ and $\lambda = 3/4$ are listed in Tables 6 and 7, respectively. In these Tables we compare the result of three methods together and with results presented by [44]. Also, the results are reported for $N = 13$ in modified generalized Laguerre functions (MGLF), $N = 20$ in Hermite functions (HF), and $N = 30$ for Sinc functions (SF). Graphs of the approximations of $f'(\eta)$ for different values of $\lambda$ MGLF, HF, and SF are shown respectively in Figures 8–10.

## 5 Conclusions

In the above discussion, we applied the collocation method to solve three nonlinear differential equations. An important concern of the collocation approach is the choice of basis functions. The basis functions have three different properties: easy computation, rapid convergence and completeness, which means that any solution can be presented to arbitrarily high accuracy by taking the truncation $N$ to be sufficiently large. We used three set of orthogonal functions as the basis function in this method and compared the results. Although all functions lead to more accurate results, it seems that based on the problem, one of the orthogonal functions goes to the more accurate result quicker.

Table 5: Comparison of $f'(0)$ for various $\lambda$ between SF collocation method and Runge-Kutta solution.

| $\lambda$ | Runge-Kutta solution [44] | SF | $h$ | $\beta$ |
|---|---|---|---|---|
| 0 | 0.94760 | 0.94749990 | 5 | 1.895 |
| $\frac{1}{4}$ | 0.91130 | 0.9100000 | 5 | 1.787 |
| $\frac{1}{3}$ | 0.90030 | 0.9003098 | 5 | 1.80062 |
| $\frac{1}{2}$ | 0.87980 | 0.8798599 | 5 | 1.75972 |
| $\frac{3}{4}$ | 0.85220 | 0.8522499 | 5 | 1.7045 |
| 1 | 0.82760 | 0.8276099 | 5 | 1.65522 |

Table 6: The comparison of $f'(\eta)$ for $\lambda = \frac{1}{4}$ for present methods and Runge-Kutta solution.

| $\eta$ | Runge-Kutta solution [44] | MGLF | HF | SF |
|---|---|---|---|---|
| 0.0 | 0.911295 | 0.911292295 | 0.911300000 | 0.9100000 |
| 0.1 | 0.813604 | 0.8136045732 | 0.818966668 | 0.8249495 |
| 0.2 | 0.721351 | 0.7214333041 | 0.739987010 | 0.7593399 |
| 0.3 | 0.635531 | 0.6357278583 | 0.671904860 | 0.6994390 |
| 0.4 | 0.556661 | 0.5570128319 | 0.612803846 | 0.6357481 |
| 0.5 | 0.484997 | 0.4854870213 | 0.561171050 | 0.5635668 |
| 0.6 | 0.420587 | 0.4211045416 | 0.515799214 | 0.4836353 |
| 0.7 | 0.363276 | 0.3636408416 | 0.475715519 | 0.4010437 |
| 0.8 | 0.312677 | 0.3127459036 | 0.440129003 | 0.3224731 |
| 0.9 | 0.268264 | 0.2679864241 | 0.408391244 | 0.2533653 |
| 1.0 | 0.229508 | 0.2288793877 | 0.379966629 | 0.1964756 |
| 1.1 | 0.195878 | 0.1949179596 | 0.354409490 | 0.1519869 |
| 1.2 | 0.166847 | 0.1655912361 | 0.331346774 | 0.1184391 |
| 1.3 | 0.141837 | 0.1403992695 | 0.310464102 | 0.0937062 |
| 1.4 | 0.120362 | 0.1188633223 | 0.291495126 | 0.0756563 |
| 1.5 | 0.102025 | 0.1005335321 | 0.274212962 | 0.0624714 |
| 2.0 | 0.043951 | 0.0436775689 | 0.207169774 | 0.0315387 |
| 2.5 | 0.018546 | 0.0196829528 | 0.162014544 | 0.0189928 |
| 3.0 | 0.007610 | 0.0088012343 | 0.130161222 | 0.0106135 |
| 3.5 | 0.002953 | 0.0030440964 | 0.106855392 | 0.0042628 |
| 4.0 | 0.000962 | $-0.0001159364$ | 0.089291515 | $-0.0006069$ |
| 4.5 | 0.000123 | $-0.0014217904$ | 0.075727337 | $-0.0043205$ |

Table 7: The comparison of $f'(\eta)$ for $\lambda = \frac{3}{4}$ for present methods and Runge-Kutta solution.

| $\eta$ | Runge-Kutta solution [44] | MGLF | HF | SF |
|---|---|---|---|---|
| 0.0 | 0.852193 | 0.852287074 | 0.852100000 | 0.8522499 |
| 0.1 | 0.755377 | 0.7553472043 | 0.760260331 | 0.7663430 |
| 0.2 | 0.665448 | 0.6652437998 | 0.682506144 | 0.6986781 |
| 0.3 | 0.582985 | 0.5826438854 | 0.616097677 | 0.6367457 |
| 0.4 | 0.508141 | 0.5077870112 | 0.558930301 | 0.5721998 |
| 0.5 | 0.440849 | 0.4406043145 | 0.509365693 | 0.5012224 |
| 0.6 | 0.380907 | 0.3808120147 | 0.466113136 | 0.4248815 |
| 0.7 | 0.327973 | 0.3279838418 | 0.428144630 | 0.3479289 |
| 0.8 | 0.281536 | 0.2816065940 | 0.394633242 | 0.2762278 |
| 0.9 | 0.241013 | 0.2411219871 | 0.364907692 | 0.2142853 |
| 1.0 | 0.205832 | 0.2059574952 | 0.338418505 | 0.1641129 |
| 1.1 | 0.175434 | 0.1755486532 | 0.314712300 | 0.1254656 |
| 1.2 | 0.149275 | 0.1493544263 | 0.293412523 | 0.0967406 |
| 1.3 | 0.126821 | 0.1268673543 | 0.274204137 | 0.0758559 |
| 1.4 | 0.107596 | 0.1076195445 | 0.256822044 | 0.0608150 |
| 1.5 | 0.091196 | 0.0911857686 | 0.241041852 | 0.0499612 |
| 2.0 | 0.039223 | 0.0392999700 | 0.180361082 | 0.0249502 |
| 2.5 | 0.016574 | 0.0165687080 | 0.140011901 | 0.0146721 |
| 3.0 | 0.006832 | 0.0067312540 | 0.111830465 | 0.0076143 |
| 3.5 | 0.002668 | 0.0026491548 | 0.091374230 | 0.0022177 |
| 4.0 | 0.000913 | 0.0011427557 | 0.076057527 | $-0.0019132$ |
| 4.5 | 0.000237 | 0.0006565451 | 0.064292470 | $-0.0050431$ |

Figure 8: MGLF approximation of $f'(\eta)$ for different values $\lambda = 0, 1/4, 1/3, 1/2, 3/4$ and 1.

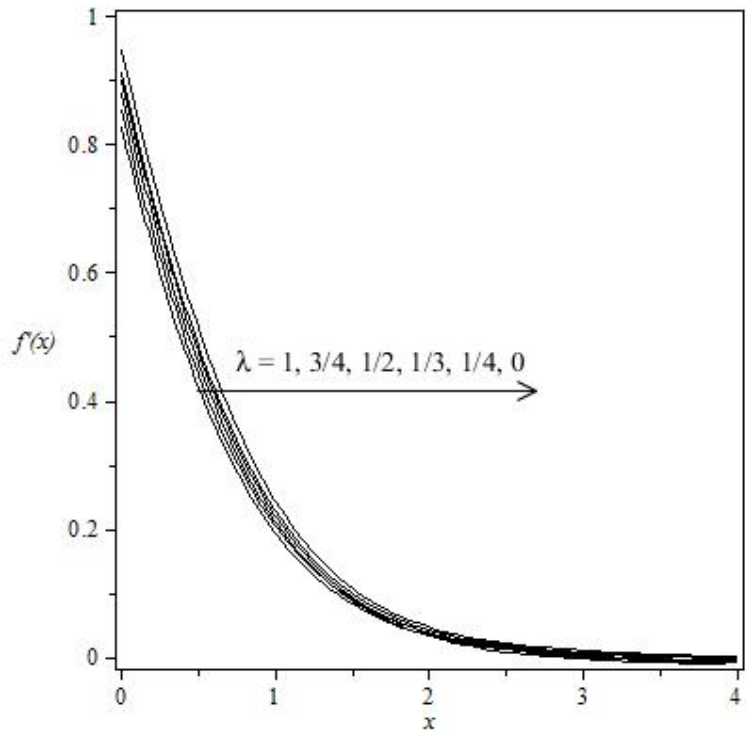

Figure 9: HF approximation of $f'(\eta)$ for different values $\lambda = 0, 1/4, 1/3, 1/2, 3/4$ and 1.

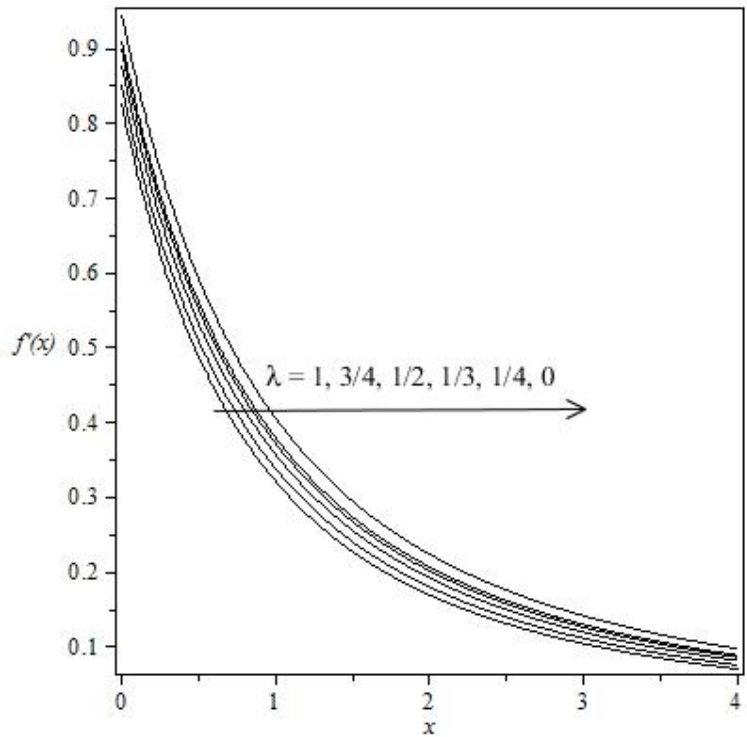

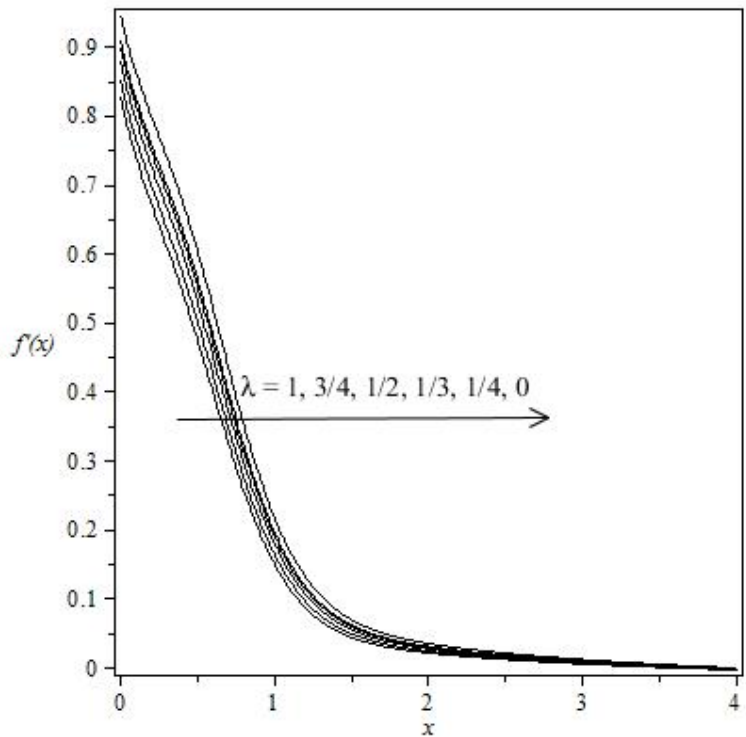

Figure 10: SF approximation of $f'(\eta)$ for different values $\lambda = 0, 1/4, 1/3, 1/2, 3/4$ and $1$.